\documentclass[11pt]{article}
\usepackage{mathpple}
\usepackage{graphicx}
\usepackage{amssymb, amsmath, latexsym}
\usepackage{picins}
\usepackage{amscd}
\newtheorem{theorem}{Theorem}

\newtheorem{definition}{Definition}
\newtheorem{lemma}{Lemma}
\newtheorem{prop}{Proposition}

\renewcommand{\mod}{\text{ mod }}

\newcommand{\proof}[1]{\textbf{Proof: }{#1}\hfill$\Box$}

\newcommand{\Out}[1]{\operatorname{Out}(F_{#1})}

\title{Ornate Necklaces and the Homology of the Genus One Mapping Class group}
\author{James Conant}
\begin{document}
\date{}
\maketitle
\begin{abstract} 
According to seminal work of Kontsevich, the unstable homology of the mapping class group of a surface can be computed via the homology of a certain lie algebra. In a recent paper, S. Morita analyzed the abelianization of this lie algebra, thereby
constructing a series of candidates for unstable classes in the homology of  the mapping class group. In the current paper, we show that these cycles are all nontrivial, representing homology classes in  $H_{k}(\mathbf M^{k}_1;\mathbb Q)_{\mathfrak S_{k}}$ for all $k\geq 5$ satisfying $k\equiv 1\mod 4$. Here $\mathbf M^k_1$ is the mapping class froup of a genus one surface with $k$ punctures. 
\end{abstract}

\section{Introduction}
Although the stable cohomology of the mapping class group has been completely computed \cite{Madsen-Weiss}, the unstable cohomology remains an interesting target of exploration. A well-known theorem of Kontsevich \cite{Kon} (see also \cite{OnKon}) relates this unstable cohomology with the homology of a certain symplectic lie algebra of derivations, $\mathfrak a_\infty$. To produce nontrivial elements in the cohomology of some lie algebra $\mathfrak a$, it can be profitable to find an abelian quotient, $\mathfrak b$, of $\mathfrak a$, which yields a map $$\bigwedge \mathfrak b=H^*(\mathfrak b)\to H^*(\mathfrak a).$$ This produces many cocycles in the image which become candidate cohomology classes. In \cite{Mor99}, Morita uses this technique to construct a sequence of cycles for the homology of $\Out{n}$, which Vogtmann and I analyzed in \cite{ConVog}, showing that the first two are nontrivial. 

The positive degree part of the lie algebra $\mathfrak a_\infty$ is the direct limit of lie algebras $\mathfrak a^+_n$. In this case, Morita \cite{Mor06} calculates the weight $2$ part of the abelianization of $\mathfrak a^+_n$, $H_1(\mathfrak a_n^+)_2$, from which he produces a sequence of cycles in the homology of the mapping class group. It is the purpose of this paper to show that these are all nonzero, and in fact represent nontrivial classes in $$H_{k}(\mathbf M^{k}_1;\mathbb Q)_{\mathfrak S_{k}}, \hspace{1em}k\geq 5, \hspace{.5em}k\equiv1\text{ mod }4$$ where $\mathbf M^{m}_g$ represents the genus $g$ mapping class group with $m$ punctures, and $\mathfrak S_m$ is the symmetric group permuting the punctures.

The homology of the genus one mapping class group has been studied previously by Getzler \cite{Getzler}. In fact, he explicitly computes the Euler characteristic of both $\mathbf M^m_1$ and $\mathbf M^m_1/\mathfrak S_m$. The generating function for the latter is
\begin{align*}
\sum_{m=1}^\infty \chi(\mathbf M^m_1)x^m&=(x+x^2+x^3)\frac{(1-x^4-2x^8-x^{12}+x^{16})}{(1-x^8)(1-x^{12})}\\
&=x+x^2+x^3-x^5-x^6-x^7-x^9-x^{10}-x^{11}-x^{13}-x^{14}-x^{15}-x^{17}-x^{18}-x^{19}\\
&-3x^{21}-3x^{22}-3x^{23}-x^{25}-x^{26}-x^{27}-3x^{29}-3x^{30}-3x^{31}-3x^{33}-3x^{34}-3x^{35}\\&-3x^{37}-3x^{38}-3x^{39}-3x^{41}-3x^{42}-3x^{43}-5x^{45}-5x^{46}-5x^{47}-3x^{49}-3x^{50}\cdots
\end{align*}
Thus the homology is growing, although it may only be growing very slowly. Thus these classes constructed from $H_1(\mathfrak a_g^+)_2$ only form a small part of the homology. The genus one mapping class group remains an interesting object of investigation!

The outline of the paper is as follows. In the first section we review the definition of ribbon graph homology, which is a useful way to compute the cohomology of the mapping class group. We then introduce cocycles, $\Theta_k$, on this complex, and show they are nontrivial classes by constructing explicit cycles, $Z_k$, with which they pair nontrivially. In the last section we explain why these cocycles correspond to the ones constructed by Morita. 

We finish the introduction with a couple of questions. In the main construction, the cocycles $\Theta_k$ are only nonzero for $k\in\{5,9,13,\ldots\}$, whereas the cycles are nonzero whenenever $k\geq 3$ is odd. 

{\bf Question:} Do the cycles $Z_k$ represent nontrivial elements of $H^k(\mathbf M^k_1;\mathbb Q)_{\mathfrak S_k}$ for all odd $k\geq 3$?

Also, these classes $\Theta_k$ only form a small part of the homology, but arise in a very elegant way via the lie algebra theory. Hence it is natural to wonder what is special about them within the homology of the mapping class group.

{\bf Question (Morita):} Can one identify these (co)cycles in the context of Getzler's \cite{Getzler} computation using algebraic geometry and number theory?

{\bf Acknowledgements:} The author wishes to thank Shigeyuki Morita for helpful discussions. 
The author was supported by NSF grant DMS-0604351.

\section{The ribbon graph complex}
We review the definition of a well-known graphical chain complex computing the homology of the mapping class group. This is the ribbon graph complex.

\begin{definition}\hspace{1em}
\begin{enumerate}
\item A \emph{ribbon graph} is a finite connected graph with vertices of valency at least $3$, with the additional structure that each vertex has a specified cyclic order of all incoming half-edges. By convention, when these graphs are drawn, they inherit these cyclic orders from the orientation of the plane of the paper. Every ribbon graph can be thickened canonically into an oriented surface.

\item An \emph{orientation} of a graph, $X$, can be defined in many equivalent ways. It is an equivalence class of decorations on a graph, and in the connected case, there are exactly two orientations. For our purposes, an orientation will be determined by an ordering of the vertices and a direction for all the edges. Two such assignments of vertex order and edge direction are equivalent iff they differ by an even number of vertex swaps and edge reversals.

\item If $m>0$, let $r\mathcal G^m_g$ be the rational vector space spanned by oriented ribbon graphs which thicken to a surface of genus $g$ with $m$ punctures, modulo the relations that
$(X,-or)=-(X,or)$, where $or$ is an orientation of the graph $X$. 

\item Let $\displaystyle r\mathcal G=\bigoplus_{g\geq0,m>0}r\mathcal G^m_g$.

\item The vector space $r\mathcal G^m_g$ is graded by number of vertices. Let the degree $k$ subspace be denoted by $r\mathcal G^m_g[k]$.

\item The boundary operator $d\colon r\mathcal G^m_g[k]\to rG^m_g[k-1] $ is defined 
by letting $d(X)$ be
the sum of oriented graphs obtained by contracting each non-loop edge of $X$.
Note that the two vertices which are joined by the edge contraction inherit a canonical cyclic order. The orientation is determined by reordering the vertices so that the contracted edge will go from vertex $1$ to vertex $2$. Then in the contracted graph, the newly created vertex is ordered first, and the the other vertices are shifted down by one from their previous numerical label.
\end{enumerate}
\end{definition}
The following result is well-known, and follows from work of Penner \cite{Penner}. See also \cite[Theorem 4]{OnKon}.
\begin{theorem}\label{PennerThm}
There is an isomorphism $$H^{4g+2m-4-k}(\mathbf M_g^m;\mathbb Q)_{\mathfrak S_m}\cong H_k(r\mathcal G_g^m).$$
\end{theorem}

It will be convenient to consider a quasi-isomorphic quotient complex $\overline{r\mathcal G}.$

\begin{definition}\hspace{1em}
\begin{enumerate}
\item A ribbon graph $X$ is said to have a \emph{cut vertex}, if in the thickening of $X$, there is a properly embedded arc which meets $X$ only at the vertex and disconnects the thickened surface into two pieces. 
\item Let $\mathcal C$ be the subcomplex spanned by graphs with cut vertices, and let $\overline{r\mathcal G}=r\mathcal G/\mathcal C$.
\end{enumerate}
\end{definition}
\begin{prop}
The natural map $r\mathcal G\to\overline{r\mathcal G}$ induces a quasi-isomorphism.
\end{prop}
\proof{This follows from remarks at the end of Section 4.2 of \cite{OnKon}. If one considers the associative operad to be spanned by planar rooted binary trees as we did in that paper, then a separating edge inside a tree at a vertex makes the vertex a cut vertex. If there is a separating edge of the graph outside of a vertex, then either of its two endpoints will correspond to cut vertices of the ribbon graph. 

One could also reproduce the algebraic argument from \cite{Cut}, where the only change needed is to slightly modify the proof of Lemma 2.1 to account for the ribbon structure.}

\section{Morita's unstable classes}
In this section, for every $k\geq 5$ satisfying $ k\equiv 1\mod 4$,
we construct a cocycle $\Theta_{k}$, on the ribbon graph complex, and show that it represents a nontrivial class in
$$H^{k}(r\mathcal G_1^{k})\cong H_{k}(\mathbf M^{k}_1,\mathbb Q)_{\mathfrak S_{k}}.$$

In the next section we will review Morita's construction and show that the classes $\Theta_k$ indeed coincide with the the ones he constructed.

\begin{figure}
\begin{center}
\includegraphics[width=.3\linewidth]{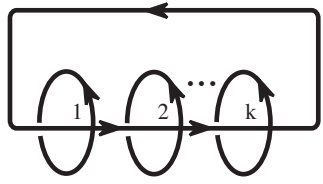}
\caption{The ribbon graph $X_k$.}\label{Xk}
\end{center}
\end{figure}

\begin{definition}
For any integer $k\geq 1$, let $X_{k}\in r\mathcal G^{k}_1$ be the graph pictured in Figure~\ref{Xk}.
\end{definition}

\begin{prop}
The ribbon graphs $X_k$, $k\geq 1$, are nonzero if and only if $k\equiv 1\mod 4$ and $k\neq 1$. 
\end{prop}
\proof{
\parpic[fr][r]{\includegraphics[width=.3\linewidth]{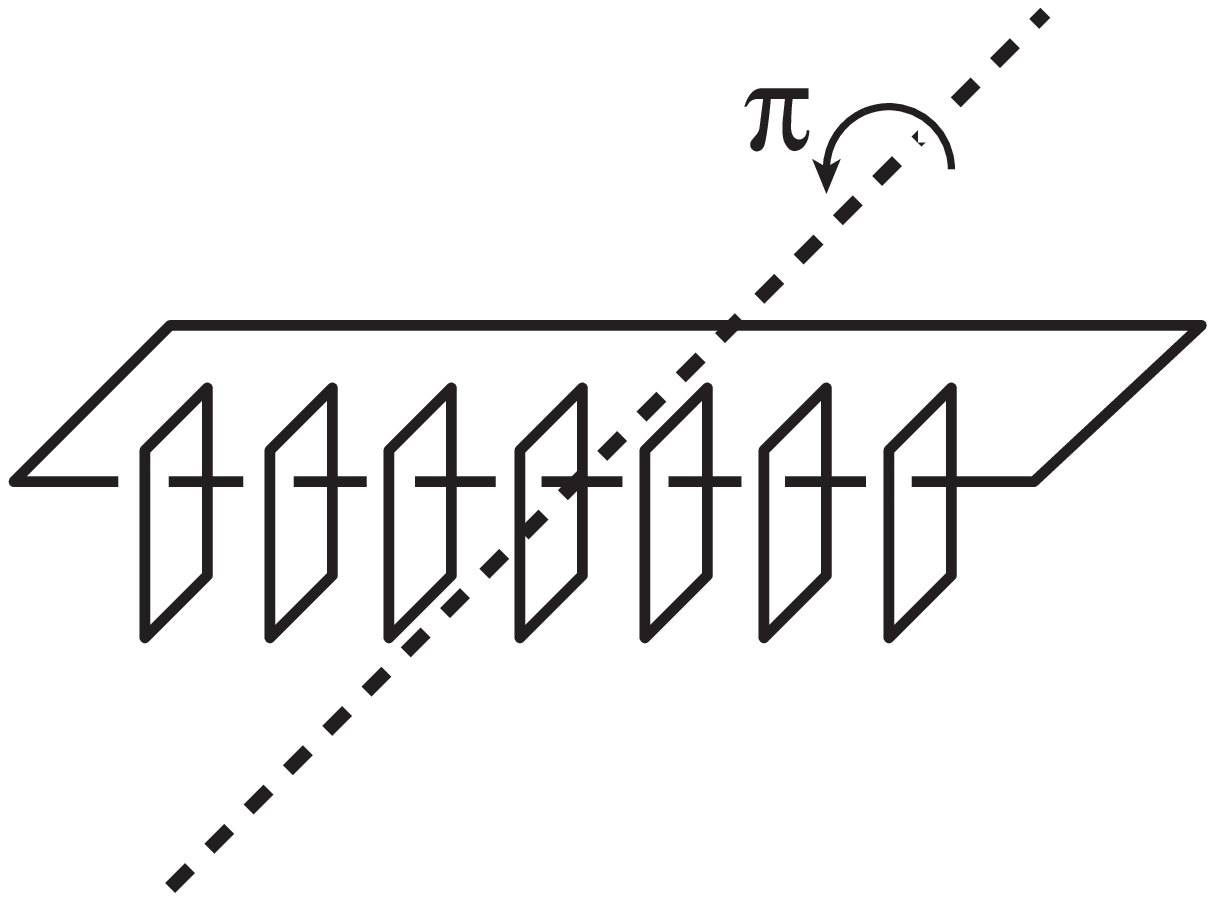}}
First I would like to argue that when $k>1$, the symmetry group of $X_k$ as an unoriented ribbon graph is the dihedral group $D_{2k}$. There is clearly a cyclic symmetry which preserves the ribbon structure. Reflection in a line is achieved by rotation of $X_k$ by $\pi$ through an axis, as in the picture on the right.
This clearly preserves the ribbon structure. Thus the dihedral group forms a subgroup of the symmetry group. Conversely any symmetry of the graph induces a symmetry of the big loop, and this gives a map from the symmetry group onto the dihedral group. 
We claim that this is a monomorphism. To see this, suppose that the induced element of the dihedral group is the identity. Then all vertices are fixed, and the only thing that could change is the order in which the ends of each small loop attach. However, these cannot change because that would change the cyclic order at the vertex, and therefore not preserve the ribbon structure. Hence the original symmetry was the identity.

In the case that $k$ is even, the generator of the cyclic symmetry group is orientation-reversing. Hence $X_k=-X_k$ and so $X_k=0$. In the case that $k\equiv 3\mod 4$, the rotation by $\pi$ around an axis, as pictured above, is orientation-reversing, and so $X_k=0$. To see this, note that the edge directions are all reversed, while the vertices are changed by an odd permutation. Since there are an even number of edges, the result follows. In the case that $X_k\equiv 1\mod 4$ and $k\neq 1$, it is easy to check that both of these types of symmetries are orientation preserving, and since they generate the entire symmetry group, the graph is nonzero.

Finally, when $k=1$, there is an orientation-reversing symmetry which exchanges the two loops. Visualize the single vertex with the edge of one loop emanating from the top and feeding into the bottom. The other loop emanates from the right and feeds into the left. If one rotates the picture by $-\pi/4$, the direction of the vertical edge gets switched. 
}

\begin{definition}If $k\equiv 1\mod 4$ and $k\neq 1$, then let $\Theta_k\colon r\mathcal G^{k}_1\to\mathbb Q$ be the characteristic function for $X_k$. 
\end{definition}

\begin{lemma}{\hspace{1em}}

\vspace{-1em}
\begin{enumerate}
\item $\Theta_k$ is a cocycle.
\item $\Theta_k$ induces a cocycle $\overline{\Theta}_k$ on $\overline{r\mathcal{G}}^k_1$.
\end{enumerate}
\end{lemma}
\proof{
Let $Y$ be a ribbon graph. We wish to show that $\Theta_k(dY)=0$. If $X_k$ appears as a summand in $dY$, then $\pm Y$ must be of the form

\centerline{\includegraphics[width=.35\linewidth]{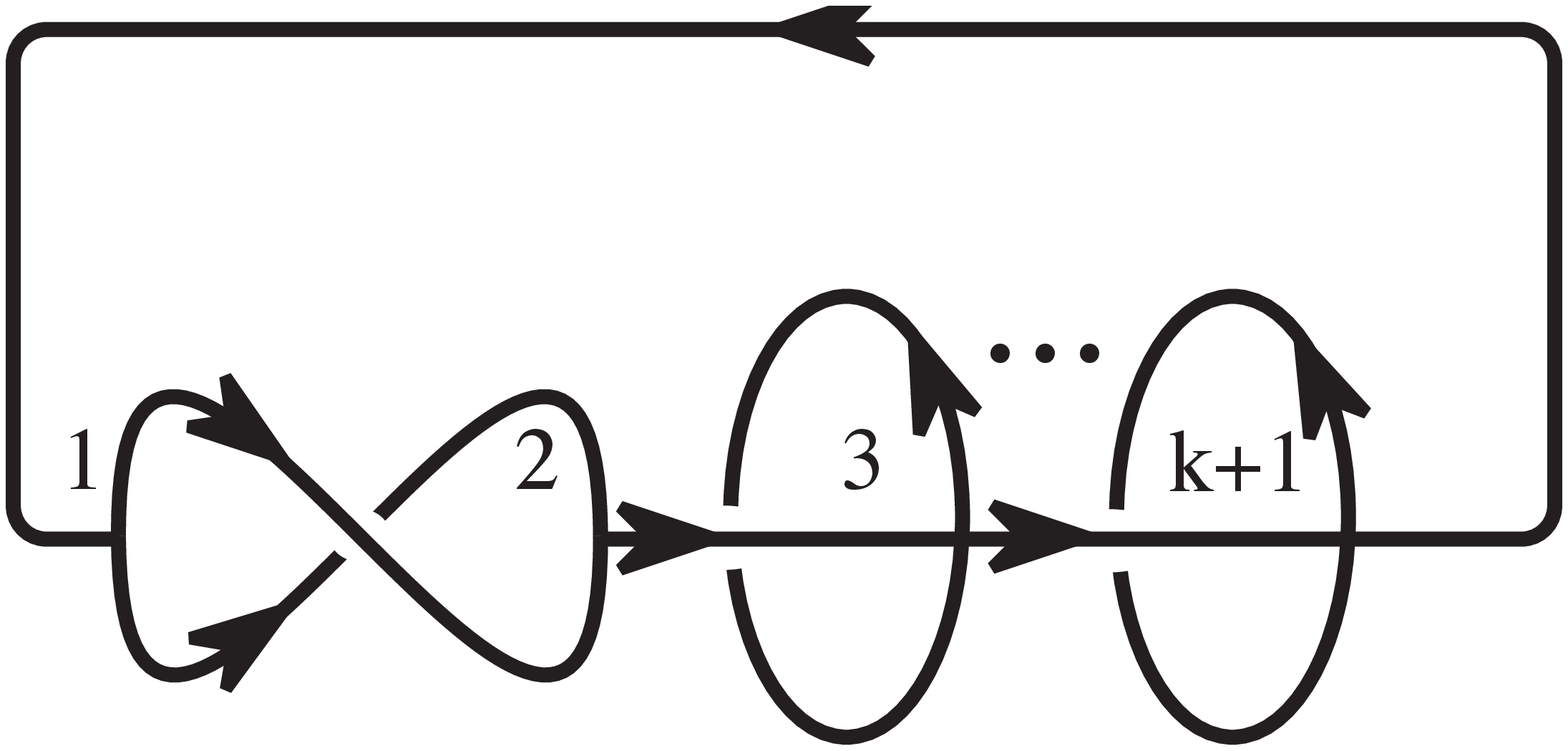}}

 But then $dY$ can be written

\centerline{\includegraphics[width=.6\linewidth]{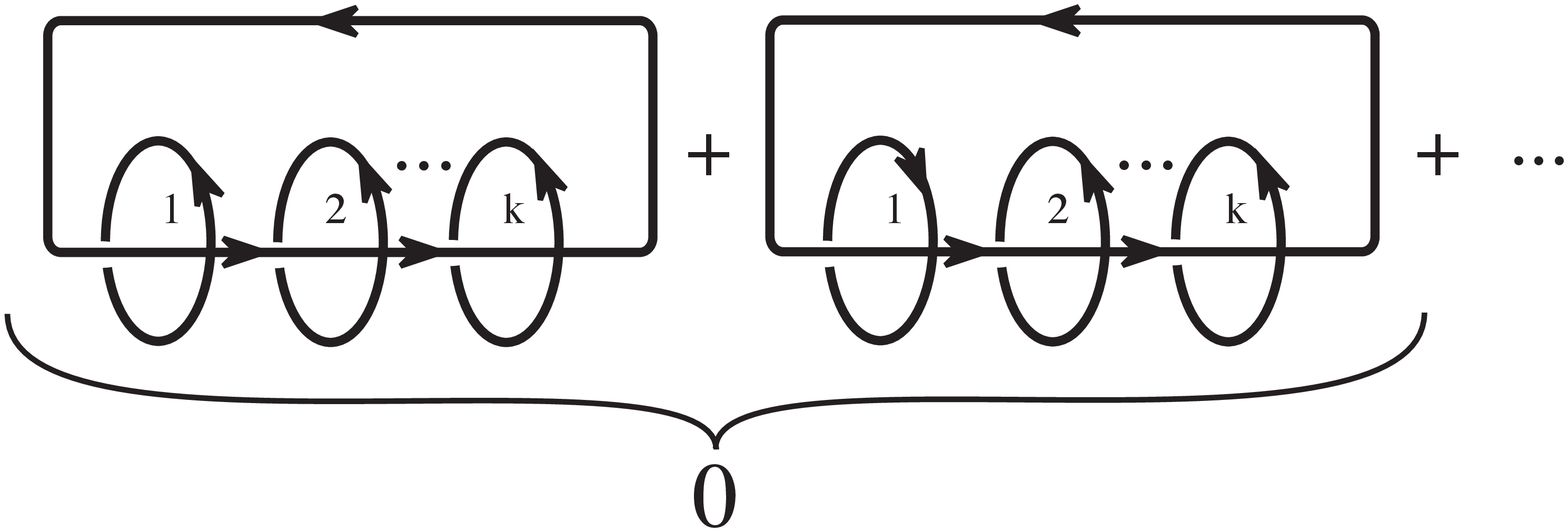}}

Thus $\Theta_k(dY)=0$ as desired.

We get a well-defined induced cocycle $\overline{\Theta}_k$ because $\Theta_k$ vanishes on all graphs with cut vertices.}

We now come to the main theorem:
\begin{theorem}
For all $k\geq 5$ where $k\equiv 1\mod 4$, $\Theta_{k}$ represents a nontrivial ribbon graph cohomology class, and therefore represents a nontrivial homology class in $$H_{k}(\mathbf M^{k}_1;\mathbb Q)_{\mathfrak S_{k}}.$$
\end{theorem}

To see this, we will find a cycle $Z_k\in\overline{r\mathcal G}^k_1$ such that $\overline{\Theta}_k(Z_k)\neq 0$.

\begin{definition}
Let $T_i$ be the sum of all isomorphism classes of planar binary rooted trees with $i$ leaves. Each such tree has a canonical orientation defined by directing the edges away from the root, and numbering the internal vertices from left to right.  
\end{definition}

\begin{figure}[ht!]

$\underbrace{\includegraphics[width=.04\linewidth]{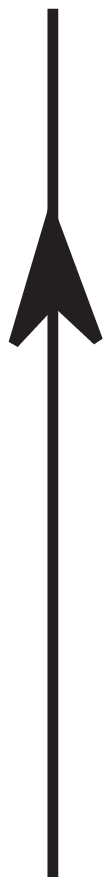}}_{T_1}$\hfill
$\underbrace{\includegraphics[width=.077\linewidth]{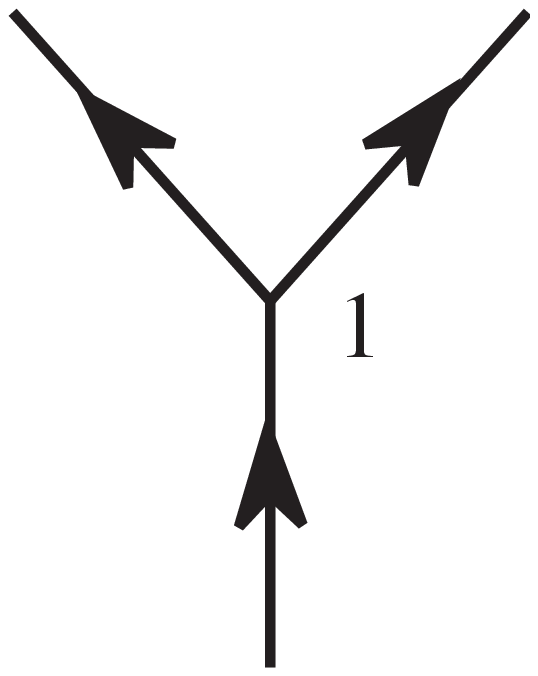}}_{T_2}$\hfill
$\underbrace{\includegraphics[width=.2\linewidth]{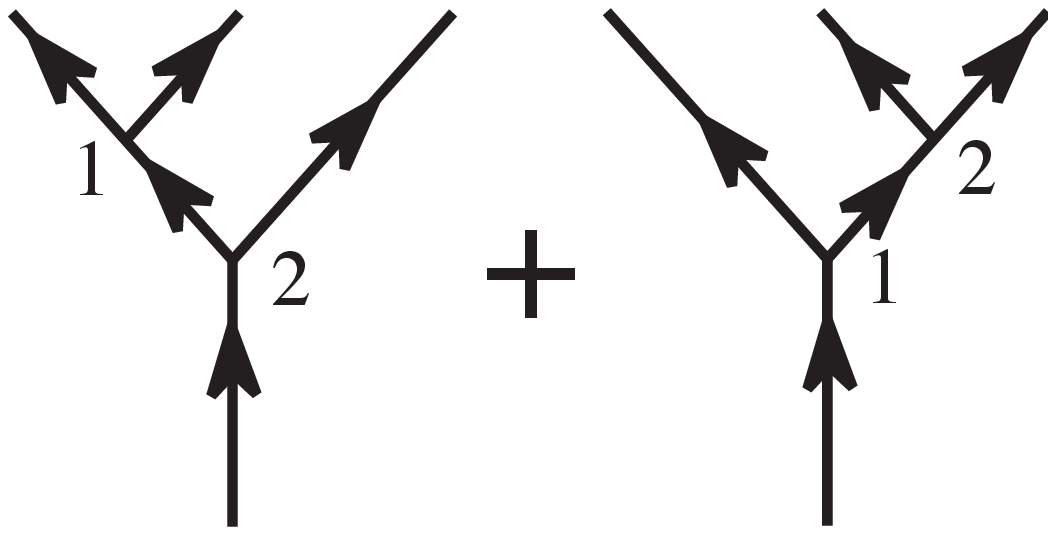}}_{T_3}$\hfill
$\underbrace{\includegraphics[width=.56\linewidth]{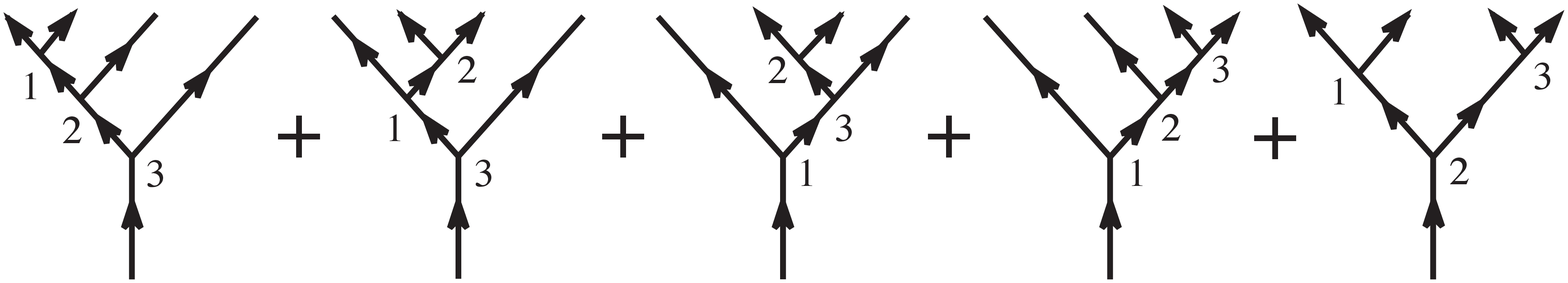}}_{T_4}$

\caption{The first few $T_i$.}\label{tifig}

\end{figure}

The first few $T_i$ are pictured in Figure~\ref{tifig}.

\begin{definition}
An \emph{ornate necklace} is a sum of ribbon graphs of the form pictured in Figure~\ref{ornate}. It is oriented so that the edges of the large loop are directed counterclockwise as indicated, and the edges in the $T_{i_j}$ are directed away from the root. The vertices are numbered so that the $T_{i_j}$ vertices lie before the $T_{i_{j+1}}$ vertices. The root of $T_{i_j}$ is numbered before its other vertices, which are ordered left to right as in the definition of $T_{i_j}$. Let the ornate necklace of Figure~\ref{ornate} be denoted $[i_1,\ldots,i_n]$.

\begin{figure}
\begin{center}
\includegraphics[width=.5\linewidth]{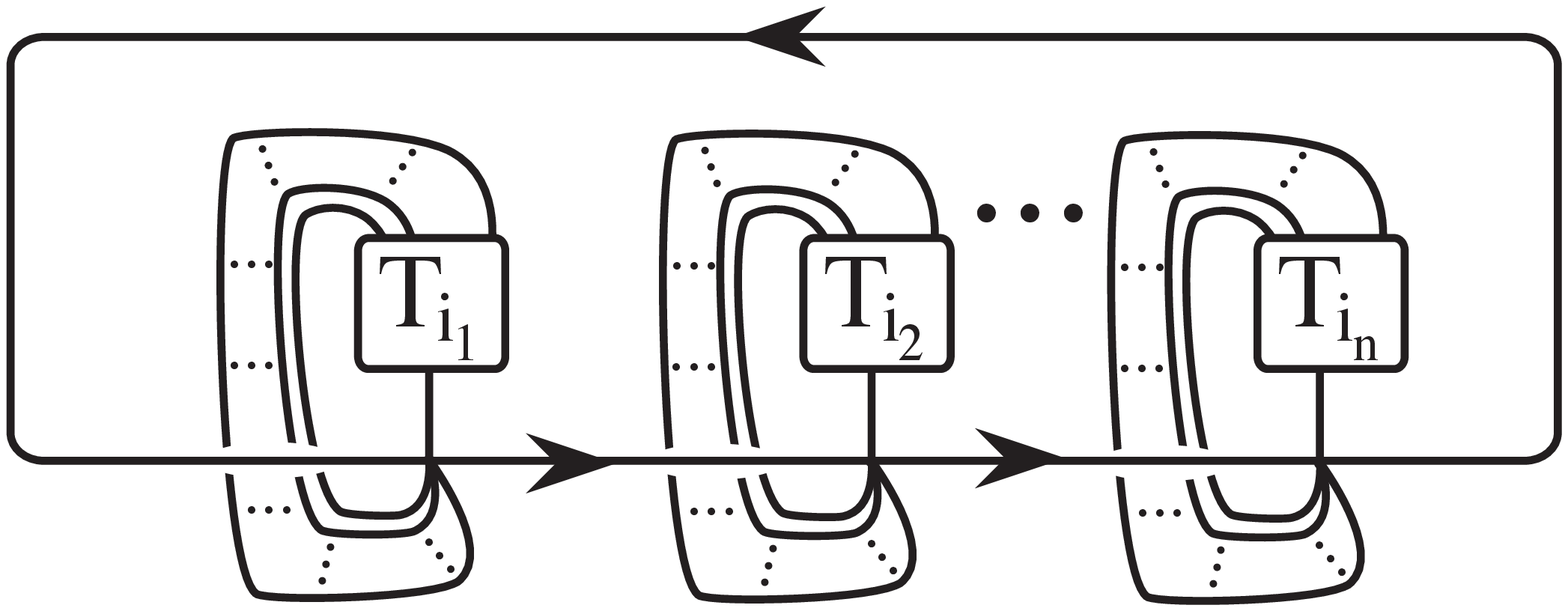}
\end{center}
\caption{The ornate necklace $[i_1,\ldots,i_n]$.}\label{ornate}
\end{figure}

\end{definition}

{\bf Examples:}
\begin{enumerate}
\item $[1,\ldots,1]=X_k$, if there are $k$ `$1$'s.
\item $[3,2]=$

\centerline{\includegraphics[width=.55\linewidth]{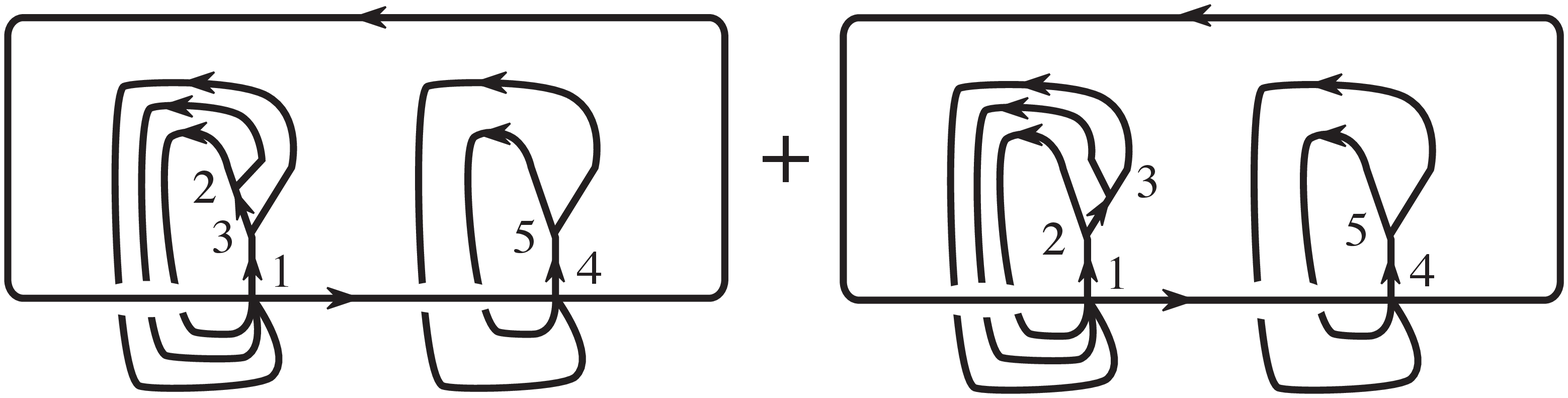}}
\end{enumerate}

\begin{definition}
The notation $(i_1,i_2)$ represents the local picture:

\centerline{\includegraphics[width=.25\linewidth]{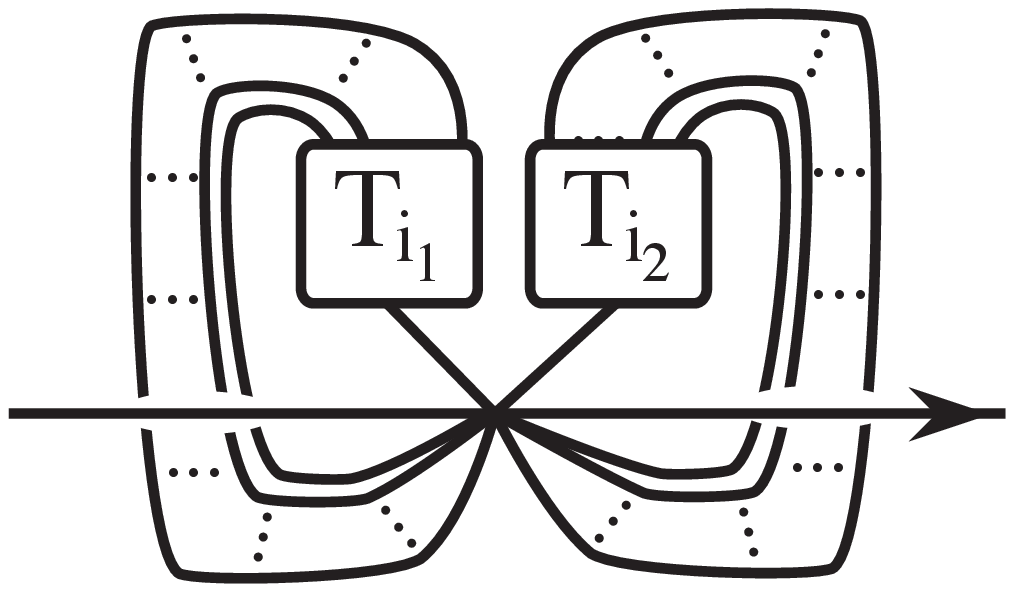}}

and can be inserted into the above bracket notation, as in
$[(i_1,i_2),i_3,\ldots,i_n]$, with evident meaning.
\end{definition}

\begin{lemma}
In $d[i_1,\ldots,i_n]$ the only edge contractions which contribute are
\begin{enumerate}
\item The root edges of each $T_{i_j}$.
\item The edges of the big loop.
\end{enumerate}
\end{lemma}
\proof{
The two types of edges not mentioned are the interior edges of the $T_{i_j}$ and the edges emanating from the tops of the $T_{i_j}$. 
\parpic[fr][r]{\includegraphics[width=.35\linewidth]{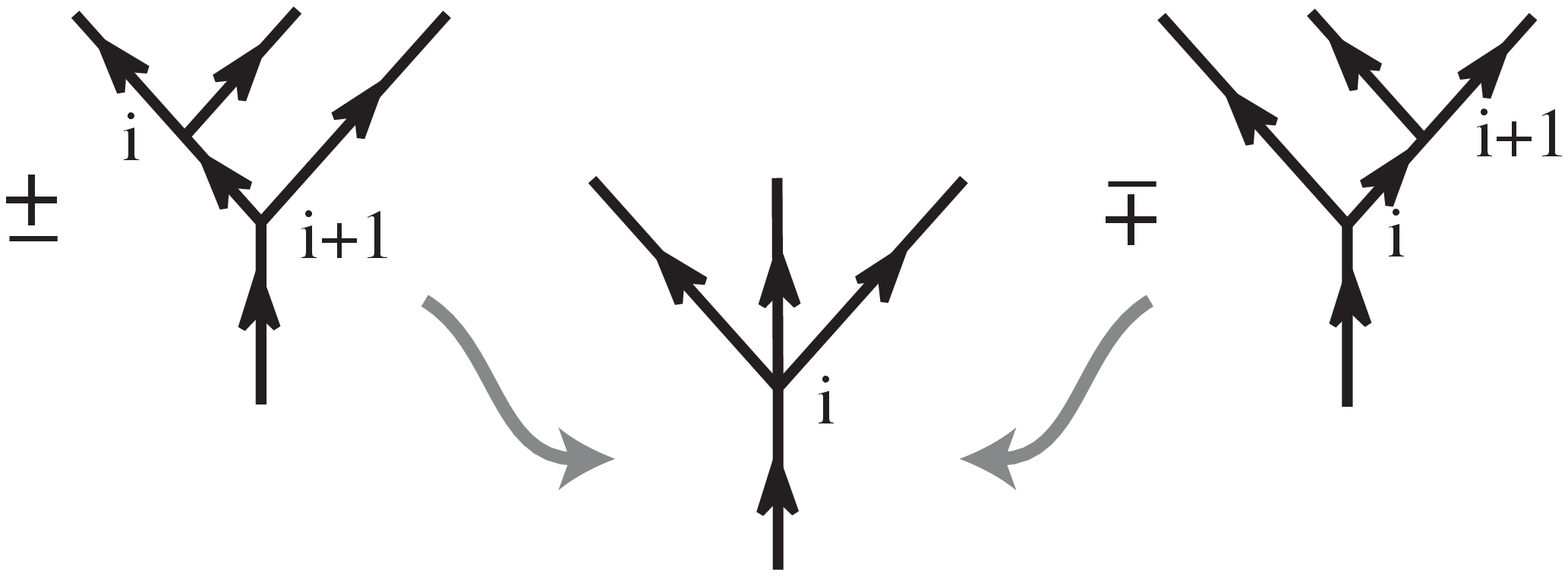}}
To see that contracting the interior edges of the $T_{i_j}$ cancel, note that such an edge contraction will create a $4$-valent vertex which can be expanded in two different ways, and thus appears twice when contracting the $T_{i_j}$ interior edges. Moreover, the signs are opposite, as indicated in the picture on the right.

For edges emanating from the top of the $T_{i_j}$, we show that the resulting ribbon graphs have cut vertices, and are therefore $0$ in $\overline{r\mathcal G}$.  Such an edge, $e$, emanates from a trivalent vertex $v$. When $e$ is contracted, the part of $T_{i_j}$ also emanating from $v$ will form a cut component at the root vertex, as indicated in the following picture.

\centerline{\includegraphics[width=.5\linewidth]{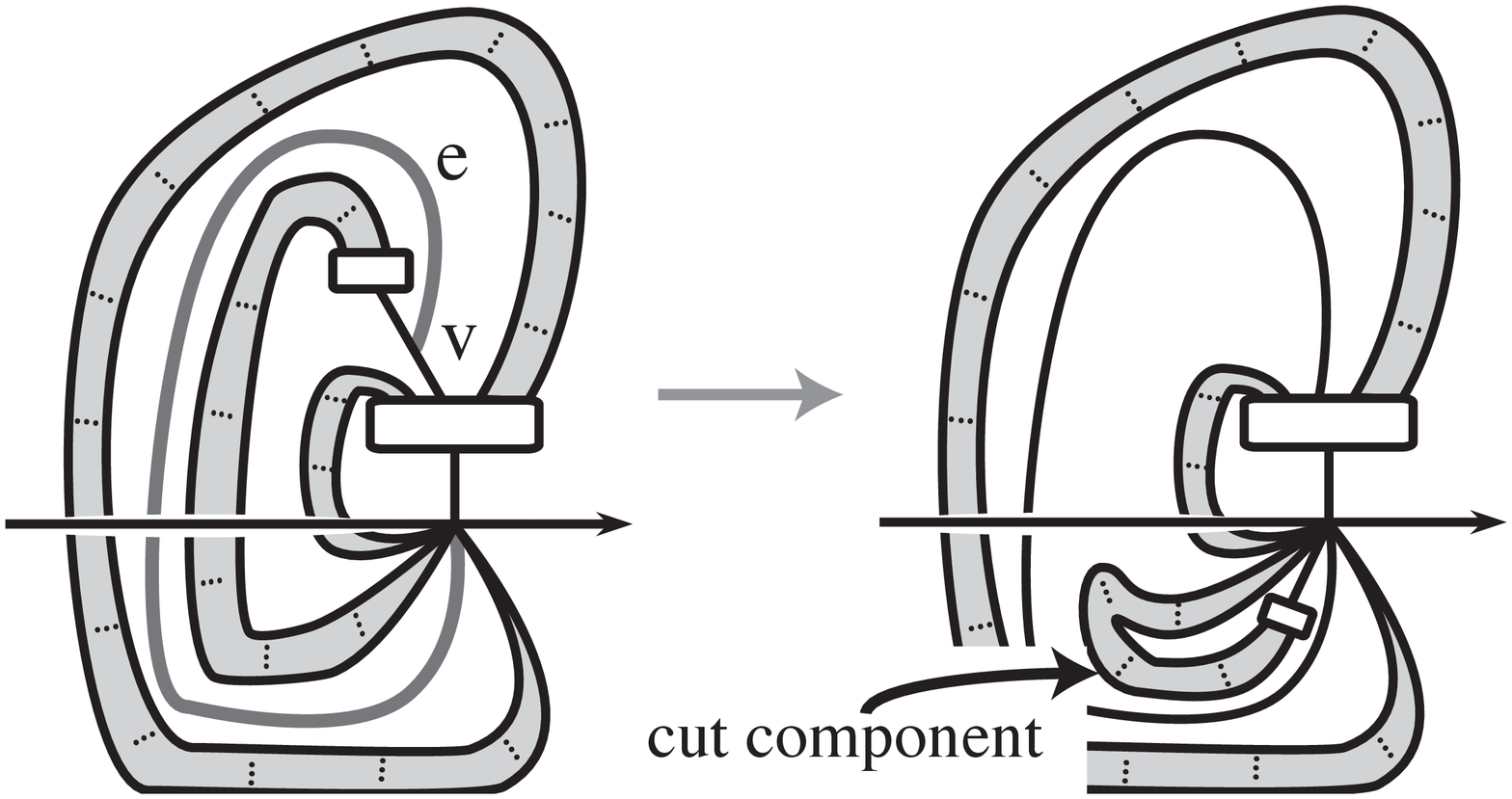}}

The grey strips represent sets of parallel edges, whereas the boxes represent trees. In the above picture, we have drawn $e$ emanating from the right of $v$, but it also evidently works if $e$ emanates from the left of $v$.
}

\begin{definition}
Let $|i_1,\ldots, i_n|$ be the number of cyclic symmetries possessed by $[i_1,\ldots,i_n]$. For example $|1,1,1,1,1|=5,|1,2,1,2,1,2|=3$ and $|1,1,2,1|=1$.
\end{definition}

Now we are ready to define the cycle $Z_k$.
\begin{definition}
Let $$Z_k=\sum (-1)^n\frac{[i_1,\ldots, i_n]}{|i_1,\ldots, i_n|},$$ where the sum is over all isomorphism classes of ornate necklaces where $i_1+\cdots +i_n=k$.
\end{definition}

This is well-defined, in the sense that the terms $[i_1,\ldots, i_n]$ have a canonical orientation invariant under cyclic symmetries. (This follows because $k$ is odd.)

For example, $$Z_5=-\frac{1}{5}[1,1,1,1,1]+[2,1,1,1]-[1,2,2]-[1,1,3]+[2,3]+[1,4]-[5].$$

Clearly, when $k\equiv 1\mod 4$ and $k\neq 1$, we have $\overline{\Theta}_k(Z_k)=-1/k$. So it suffices to show $d(Z_k)=0$.

\begin{prop}
$d(Z_k)=0$.
\end{prop}
\proof{
\parpic[fr][r]{\includegraphics[width=.38\linewidth]{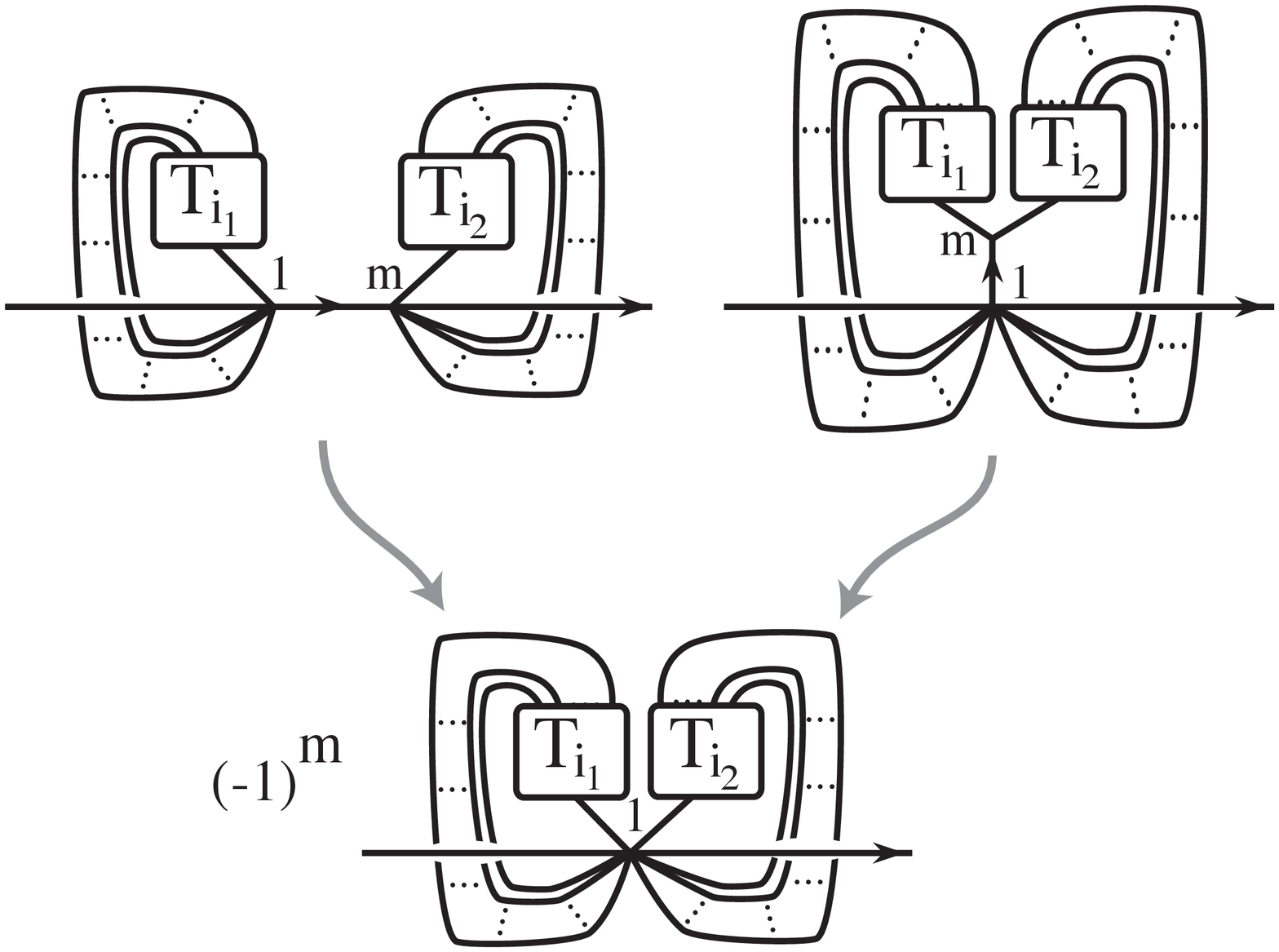}}
By the previous lemma, $d(Z_k)$ has two types of terms. 
 If  an edge in the large loop is contracted, two adjacent $T_{i_j}$'s will be joined:
$$[\ldots,i_j,i_{j+1},\ldots ]\to[\ldots,(i_j,i_{j+1}),\ldots].$$
Note that we can apply a cyclic symmetry to ensure that it of this form, and we don't need to join the last tree to the first tree. This simplifies sign considerations and notation. 
If you contract a root edge, a tree $T_{i_j}$ will become $$\pm\sum_{a+b=i_j}(a,b).$$

Thus, to show $d(Z_k)=0$, it suffices to show that each isomorphism class $\alpha=[(i_1,i_2),i_3,\ldots,i_n]$ appears twice with opposite sign. Well, $\alpha$ arises from two terms $A=[i_1,i_2,i_3,\ldots,i_n]$ and $B=[i_1+i_2,i_3,\ldots, i_n]$. Notice that $\alpha$ appears once in $(-1)^n\frac{dA}{|A|}$ and in $(-1)^{n-1}\frac{dB}{|B|}$. Thus it suffices to note that it appears with the same sign in $dA$ as it does in $dB$. This is illustrated in the picture in the above right. Note that $m$ is the same in both cases: the number of vertices of $T_{i_1}$ plus one.}

\section{Comparison to Morita's definition}
We start by defining the lie algebra $\mathfrak a_n$ for any positive integer $n$. 

\begin{definition}
Let $V_n$ be the symplectic vector space which is the rational homology of a genus $n$ surface. Explicitly, let $V_n$ be the $2n$ dimensional vector space with basis
$p_1,\ldots, p_n, q_1,\ldots, q_n$ and symplectic form $\omega$ defined so
that $\omega (p_i,q_i)=1=-\omega(q_i,p_i)$ and is trivial on all other pairs of basis vectors. 
\end{definition}

One way to define $\mathfrak a_n$ is given in the following way \cite[Prop. 2]{Mor06}.

\begin{definition} Let $$(\mathfrak a_n)_i=(V_n^{\otimes i+2})^{\mathbb{Z}_{i+2}},$$
where $\mathbb Z_{i+2}$ acts by cyclic permutation,
 and let $\mathfrak a_n=\oplus_{i=1}^\infty (\mathfrak a_n)_i$. This defines $\mathfrak a_n$ as a vector space. The bracket $[\cdot,\cdot]\colon (\mathfrak a_n)_i\otimes(\mathfrak a_n)_j\to(\mathfrak a_n)_{i+j}$ is induced by the contraction $C_{1,i+3}\colon V_n^{\otimes i+2}\otimes V_n^{\otimes j+2}\to V_n^{\otimes i+j+2}$. (Recall that the contraction $C_{i,j}\colon V^{\otimes \ell}\to V^{\otimes \ell-2}$ is defined by
 $C_{i,j}(v_1\otimes\cdots\otimes v_\ell)=\omega(v_i,v_j)v_1\otimes\cdots \otimes\widehat{v_i}
 \otimes\cdots\otimes\widehat{v_j}\otimes \cdots\otimes v_\ell.$)
\end{definition}

We note that $\mathfrak a_n$ is isomorphic to the lie algebra $\mathcal{LO}_n$ associated to the associative operad, defined in \cite[Section 2.4.1]{OnKon}. The difference is that in \cite{OnKon}, we considered the space of coinvariants $(V_n^{\otimes i+2})_{\mathbb{Z}_{i+2}}$ instead of the isomorphic space of invariants. 

\begin{definition}\hspace{1em}
\begin{enumerate}
\item Let $\mathfrak a^+_n=\oplus_{i=2}^\infty (\mathfrak a_n)_i$.
\item Let $\displaystyle \mathfrak a^+_\infty =\lim \{\mathfrak a_1^+ \subset \mathfrak a_2^+\subset \cdots\}.$
\end{enumerate}
\end{definition}

The cohomology $H^*(\mathfrak a^+_\infty)$ has a Hopf algebra structure, as described in \cite[Prop. 7]{OnKon}, so that we can consider the primitive elements, denoted with the prefix "P". We also have natural maps $PH^k(\mathfrak a^+_{n+1})^{\mathfrak{sp}}\to
PH^k(\mathfrak a^+_{n})^{\mathfrak{sp}},$ allowing us to consider the limit $PH^k(\mathfrak a^+_\infty)^{\mathfrak{sp}}$.

\begin{definition}\hspace{1em}
\begin{enumerate}
\item Let $\overline{\mathcal O\mathcal G}$ be the reduced associative graph complex, as described in \cite{OnKon}. (Recall the adjective ``reduced" means all vertices have valency at least $3$.)
Note that $r\mathcal G\subset \overline{\mathcal{OG}}$ is the subspace of connected graphs.
\item For $t\geq 0$,
 let $\overline{\mathcal O\mathcal G}_{2t}$ be the subcomplex of graphs such that the number of edges minus the number of vertices is equal to $t$. This is the \emph{weight} $2t$ part of the complex, which induces the weight $2t$ part of the homology or cohomology, also indicated with a subscript. Note that $\overline{\mathcal O\mathcal G}_{2t}$ is finite dimensional. 
\item Let
$$\psi_n\colon\bigwedge\mathfrak a^+_n\to \overline{\mathcal{OG}}$$
be defined as in \cite[Section 2.5.2]{OnKon}, using the identification of $\mathcal{LO}_n$ with $\mathfrak{a}_n$, and let $\psi_\infty$ be the limit map. 
\item For every $t\geq 0$, let the induced map on the weight $2t$ part be denoted $(\psi_\infty)_{2t}\colon  \overline{\mathcal{OG}}_{2t}\to
\left(\bigwedge^v\mathfrak a^+_\infty\right)_{2t}$. (The fact that the image lies in the weight $2t$ part is easily checked.)
\end{enumerate}
\end{definition}

\begin{theorem}
For every $t\geq 0,$ the dual map $(\psi_\infty)^*_{2t}$ induces an isomorphism $PH^k(\overline{\mathcal{OG}})_{2t}\to PH^k(\mathfrak a^+_\infty)^{\mathfrak{sp}}_{2t}.$ 
\end{theorem}
\proof{This follows from \cite[Corollary 5]{OnKon}, with some modification. A significant difference is that we are considering 
 the lie algebra $\mathfrak a^+_\infty$ and not  the lie algebra $\mathfrak a_\infty$ as in \cite{OnKon}. This has the effect of eliminating bivalent vertices from $\mathcal{OG}$, allowing us to consider $\overline{\mathcal{OG}}$ instead, but since
\cite[Proposition 8]{OnKon} is no longer true one must consider the space of $\mathfrak{sp}$-invariants
$PH^k(\mathfrak a^+_\infty)^{\mathfrak{sp}}$ instead of simply $PH^k(\mathfrak a^+_\infty)$.
 Finally, to get the exact statement above, one restricts to the weight $2t$ part and takes the dual.
}

Since we have $$PH^k(\overline{\mathcal{OG}})\cong H^k(P\overline{\mathcal{OG}})\cong H^k(r\mathcal G),$$
Theorem \ref{PennerThm}
implies the following suitably modified result of Kontsevich \cite{Kon}.

\begin{theorem}
We have $$ PH^k(\mathfrak a^+_\infty)^{\mathfrak{sp}}_{2t}\cong\bigoplus_{\scriptsize\begin{array}{c}2g-2+m=t\\ m>0\end{array}}H_{2t-k}(\mathbf M^m_g;\mathbb Q)_{\mathfrak S_m}.$$
\end{theorem}

Morita precisely determines $H_1(\mathfrak a^+_n)_2$ as follows \cite[Theorem 6]{Mor06}.
\begin{definition}
 Let $\mathfrak b_n=\bigwedge^2 H/\mathbb Q(\omega_0)$, where $\omega_0$ is the symplectic element. Let $\mathfrak b_n$ be given an abelian lie algebra structure.
\end{definition} 
 \begin{prop}[Morita \cite{Mor06},Thm. 6] 
 There is an isomorphism of $\mathfrak{sp}(2n)$ modules,
$H_1(\mathfrak a^+_n)_2\cong \mathfrak b_n$, which is induced by the map $h_1\otimes h_2\otimes h_3\otimes h_4\mapsto \omega(h_1,h_3) h_2\wedge h_4$. 
\end{prop}
\begin{definition}\hspace{1em}
\begin{enumerate}
\item Let $\pi_n\colon\mathfrak a^+_n\to \mathfrak b_n$ be the corresponding map of lie algebras.
\item There are maps $\mathfrak{b}_n\to\mathfrak{b}_{n+1}$ defined via the maps
$H_1(\mathfrak a^+_n)_2\to H_1(\mathfrak a^+_{n+1})_2$. Let $\mathfrak b_\infty$ be the direct limit, with limit map $\pi_\infty\colon \mathfrak a^+_\infty\to \mathfrak b_\infty$.
The cohomology of $\mathfrak b_\infty$ also forms a
Hopf algebra. 
\item Suppose $k\in\{5,9,13,\ldots\}$.
Let $\xi_{k,n}\colon\bigwedge^k(\mathfrak b_n)\to\mathbb Q$ be the map induced by the product of contractions: $$C_{2,3}C_{4,5}\cdots C_{2k-1,2k}C_{2k,1}\colon V_n^{\otimes 2k}\to \mathbb Q,$$ where $\bigwedge^2V_n$ is regarded as a subspace of $V_n^{\otimes 2}$ and
$\bigwedge^k(\bigwedge^2 V_n)$ is then regarded as a subspace of $V_n^{\otimes 2k}$.

Let $\xi_k\colon \bigwedge^k(\mathfrak b_\infty)\to\mathbb Q$ be the limit map. One must check that $\xi_{k,n}$ is well-defined, which amounts to showing that $\xi_{k,n}$ vanishes on $\omega_0\wedge x_2\wedge\cdots\wedge x_k$, which is easy to check provided $k\neq 1$.
\end{enumerate}
\end{definition}

\begin{prop}[Morita \cite{Mor06},Prop. 10]
We have that $$PH^k(\mathfrak b_\infty)^{\mathfrak{sp}}=\begin{cases}\mathbb Q& k=5,9,13,\ldots\\
0&\text{ otherwise}\end{cases}$$
Moreover, for $k=5,9,13.\ldots$, the map $\xi_k\colon P\bigwedge^k\mathfrak b_\infty\to\mathbb Q$
is a generator.
\end{prop}
{\bf Sketch of proof:}
One calculates 
$P\bigwedge^k(\mathfrak b_n)^{\mathfrak{sp}}=\begin{cases}\mathbb Q& k=5,9,13,\ldots\\
0&\text{ otherwise}\end{cases}$ by classical invariant theory. The single class, $e^k_n$, in degree $k$ comes from a $k$-gon and is defined as $$e^k_n
=\sum_{i=1}^n[p_i\wedge q_i]\wedge\cdots\wedge[p_i\wedge q_i],$$
where $[p_i\wedge q_i]$ is the element in $\mathfrak b_n$ represented by $p_i\wedge q_i\in \bigwedge^2 V_n$.

To see that $\xi_k$ is a nonzero invariant, we show that each $\xi_{k,n}$ has this property.
It is straightforward to show that $\xi_{k,n}$ is killed by the $\mathfrak{sp}$ action, and it is nonzero, for example, since $\xi_{k,n}(e^k_n)=-2n$.
\hfill$\Box$ 

\begin{definition}
Morita's classes are defined from the weight $2k$ classes $\xi_k$, as follows.
Note that the map $\pi_\infty$ induces   $(\pi_\infty)^*_{2k}\colon PH^k(\mathfrak b^\infty)_{2k}^{\mathfrak{sp}}\to PH^k(\mathfrak a^+_\infty)_{2k}^{\mathfrak{sp}}$. The cocyles constructed by Morita are then
$(\pi_\infty)_{2k}^*(\xi_k)$.
\end{definition}

We are now ready to relate Morita's cocycles with the cocycles $\Theta_k$ defined earlier in the paper.

\begin{prop}
We have $(\pi_\infty)^*_{2k}(\xi_k)=\pm(\psi_\infty)^*_{2k}(\Theta_k).$
\end{prop}
\proof{
Let $x\in(\bigwedge^k\mathfrak a^+_\infty)_{2k}$, be a wedge of $k$ symplecto-spiders of weight $2k$, using the terminology of \cite{OnKon}. 
Then
$$\langle (\psi_\infty)^*_{2k}(\Theta_k),x\rangle=\langle(\pi_\infty)^*_{2k}(\xi_k),x\rangle$$ if and only if $$\langle \Theta_k,\psi_\infty(x)\rangle=\langle \xi_k,\pi_\infty(x)\rangle.$$
Now $\psi_\infty(x)$ glues the labeled edges of $x$ together, multiplying by the contraction of the coefficients, and $\Theta_k$ will only be nonzero if the gluing matches opposite edges on each of the degree $2$ symplecto-spiders, and then joins the symplecto-spiders up in a $k$-cycle.
On the other hand $\pi_\infty$ by definition contracts opposite edges in each symplecto-spider, and then $\xi_k$ contracts the remaining edges into a $k$-cycle in all possible ways.
}

\end{document}